\scrollmode 
\documentclass[10pt]{llncs}
\usepackage{amssymb}
\usepackage{graphicx}
\usepackage{psfrag}
\newcommand{\mathscr }{\mathcal}
\newcommand{\vv}{{\mathbf {v}}}

\newcommand{\uu}{{\mathbf {u}}}

\newcommand{\R}{{\mathbb R}}

\newcommand{\C}{{\mathbb C}}
\newcommand{\Z}{{\mathbb Z}}
\newcommand{\Q}{{\mathbb Q}}

\newcommand{\V}{{\mathscr V}}
\newcommand{\BBB}{{\mathscr B}}
\newcommand{\WWW}{{\mathscr W}}

\newcommand{\MMM}{{\mathscr M}}

\renewcommand{\O}{{\mathscr O}}
\newcommand{\Proj}{{\mathbb P}}

\newcommand{\hp}{{\mathfrak H}}

\newcommand{\eqref}[1]{{(\ref{#1})}}
\newcommand{\Sign}{{\mathrm{sgn}}}

\newcommand{\SL}{{\mathrm{SL}}}
\newcommand{\GL}{{\mathrm{GL}}}
\newcommand{\Sp}{{\mathrm{Sp}}}
\newcommand{\SO}{{\mathrm{SO}}}
\newcommand{\supp}{{\mathrm{supp\ }}}
\newcommand{\Span}{{\mathrm{span\ }}}
\newcommand{\Gal}{{\mathrm{Gal}}}
\newcommand{\Vor}{Vorono\v{\i}}
\newcommand{\cones}{\mathbf{C}}
\renewcommand{\tilde}{\widetilde}
\begin{document}
\title{Modular symbols and Hecke operators}
\titlerunning{Modular symbols and Hecke operators}  
%
\author{Paul E. Gunnells}
\authorrunning{Paul E. Gunnells}   
%
\tocauthor{Paul E. Gunnells (Columbia University)}
\institute{Columbia University, New York, NY 10027, USA}

\maketitle              

\begin{abstract}
We survey techniques to compute the action of the Hecke operators
on the cohomology of arithmetic groups.  These techniques can be
seen as generalizations in different directions of the classical
modular symbol algorithm, due to Manin and Ash-Rudolph.  Most of the
work is contained in papers of the author and the author with Mark
McConnell.  Some results are unpublished work of Mark McConnell and
Robert MacPherson.  
\end{abstract}
\section{Introduction}
\subsection{}
Let $G$ be a semisimple algebraic group defined over $\Q $, and let
$\Gamma \subset G (\Q )$ be an arithmetic subgroup.  The cohomology of
$\Gamma $ plays an important role in number theory, through its
connection with automorphic forms and representations of the absolute
Galois group $\Gal (\bar \Q /\Q )$.  This relationship is revealed in
part through the action of the \emph{Hecke operators} on the complex
cohomology $H^{*} (\Gamma ;\C)$.  These are endomorphisms induced from
a family of correspondences associated to the pair $(\Gamma ,G (\Q
))$; the arithmetic nature of the cohomology is contained in the
eigenvalues of these linear maps.

For $\Gamma \subset \SL_{n} (\Z )$, the \emph{modular symbols} and
\emph{modular symbol algorithm} of Manin \cite{manin} and Ash-Rudolph
\cite{ash.rudolph} provide a concrete method to compute the Hecke
eigenvalues in $H^{\nu } (\Gamma;\C )$, where $\nu = n (n+1)/2-1$ is
the top degree (\S\ref{class.mod.sym}).  These symbols have allowed
many researchers to fruitfully explore the number-theoretic
significance of this cohomology group, especially for $n=2$ and $3$
\cite{agg,apt,exp.ind,geemen.top,vgt2}.  For all their power, though,
modular symbols have limitations: 
\begin{itemize}
\item The group $G$ must be the linear group $\SL _{n}$.
\item The cohomology must be in the top degree $\nu $.
\item The group $\Gamma $ must be a subgroup of $\SL _{n} (\Z )$, or
more generally  $\SL _{n} (R)$, where $R$
is a euclidean ring of integers of a number field.
\end{itemize}

\subsection{}
In this article we discuss new techniques to compute the Hecke action
on the cohomology of arithmetic groups that can be seen as
generalizing the modular symbol algorithm by relaxing the three
restrictions above.  First in \S\ref{sympms.section} we relax the
first restriction of the by replacing the linear group $\SL_{n}$ with
the symplectic group $\Sp_{2n}$ \cite{sympms}.  Next in
\S\ref{below.coho}, we relax the second restriction and consider
computations in $H^{\nu -1} (\Gamma )$, where $\Gamma \subset \SL _{n}
(\Z )$ and $n\leq 4$ \cite{experimental}.  Finally, in the last two
sections we relax all three restrictions, and consider arithmetic
groups associated to \emph{self-adjoint homogeneous cones}
(\S\ref{sahc.section}) \cite{msa,sahc}, and arithmetic groups for
which a \emph{well-rounded retract} is defined (\S\ref{wrr.section})
\cite{bobnmark}.  The first class includes $\SL_{n} (\O _{K})$, where
$\O _{K}$ is the maximal order of a totally real or CM field, as well
as arithmetic groups associated to the positive-definite $3\times 3$
Hermitian octavic matrices.  The second class includes arithmetic
subgroups of $\SL_{n} (D)$, where $D$ is a division algebra over $\Q
$.

Most of this work is contained in papers of the author
\cite{sympms,msa,experimental} or the author in joint work with Mark
McConnell \cite{sahc}.  The last section is a summary of unpublished
results of Robert MacPherson and Mark McConnell \cite{bobnmark}.  We
have omitted other work, notably that of Bygott \cite{bygott},
Teitelbaum \cite{teit}, and Merel \cite{merel}, because of lack of
space and/or author's expertise.  It is a pleasure to thank Avner Ash,
Robert MacPherson, and Mark McConnell for many conversations about
these topics.

\section{Classical modular symbols}\label{class.mod.sym}
\subsection{}\label{description}
We begin by recalling the classical modular symbol algorithm following
Ash-Rudolph \cite{ash.rudolph}.  For simplicity we consider subgroups
of $\SL_{n} (\Z )$, although everything we say can be generalized to
subgroups of $\SL_{n} (R)$, where $R$ is a euclidean maximal order in
a number field (cf. \cite{crem}).

Let $\Gamma \subset \SL_{n} (\Z )$ be a torsion-free finite-index
subgroup, and let $m\in M_{n} (\Q )$, the $n\times n$ matrices over
$\Q $.  We want to show how to use $m$ to construct a class in $H^{\nu
} (\Gamma)$.  To this end, let $X$ be the symmetric space $\SL_{n} (\R
)/\SO (n)$, let $\bar X$ be the bordification constructed by
Borel-Serre \cite{borel.serre}, and let $\partial \bar X = \bar X
\smallsetminus X$.  Let $M = \Gamma \backslash X$, $\bar M = \Gamma
\backslash \bar X$, and $\partial \bar M = \bar M \smallsetminus M$.
Then $\bar M$ is a smooth manifold with corners, and $H^{*} (\Gamma )
\cong H^{*} (\bar M)$.  We have an exact sequence
\begin{equation}\label{ex.seq}
H_{n-1}(\partial \bar X)\rightarrow H_{n}(\bar X,\partial \bar
X)\rightarrow H_{n}(\bar M,\partial \bar M)\rightarrow H^{\nu }(\bar M)
\end{equation}
coming from the sequence of the pair $(\partial \bar X, \bar X)$, the
canonical projection $\bar X \rightarrow \bar M$, and Lefschetz
duality.  Moreover, the boundary $\partial \bar X$ has the homotopy
type of the \emph{Tits building} $\BBB = \BBB _{\SL}$ associated to
$\SL_{n} (\Q)$.  This is an $(n-1)$-dimensional simplicial complex
whose $k$-simplices $\Delta $ are in bijection with flags $F$ of
rational subspaces
\[
F = \{0 \subsetneq F_{1} \subsetneq \cdots  \subsetneq
F_{k+1}\subsetneq \Q ^{n}\}; 
\]
we have $\Delta\subset \Delta'$ if and only if $F\subset
F'$.

Any ordered tuple of nonzero rational vectors determines a maximal
rational flag by defining $F_{k}$ to be the span of the first $k$ vectors.
Hence if $m\in M_{n} (\Q )$ has nonzero columns, the different
orderings of the columns determine $n!$ different oriented
$(n-1)$-simplices in $\BBB $.  These simplices can be thought of as an
oriented simplicial cycle giving a class $[m]\in H_{n-1} (\BBB) \cong
H_{n-1} (\partial\bar X)$.  The class $[m]$ is called a \emph{modular
symbol}, and these classes span $H_{n-1} (\BBB)$.  According to
Ash-Rudolph, the map $\Phi \colon H_{n-1} (\BBB) \rightarrow H^{\nu }
(\Gamma )$ induced by \eqref{ex.seq} is surjective; hence the (duals
of) the modular symbols span $H^{\nu } (\Gamma )$.

\subsection{}\label{rels}
Write $[m] = [m_{1},\dots ,m_{n}]$, where each column $m_{i}\in \Q
^{n}\smallsetminus \{0 \}$, and let $\MMM _{n}$ be the $\Z $-module
generated by the classes of the symbols $[m]$.  Using the description
in \S\ref{description}, one can show that elements of $\MMM _{n}$
satisfy the following relations:
\begin{enumerate}
\item $[qm_{1},m_{2},\dots ,m_{n}] = [m]$, for $q\in \Q ^{\times }$.
\item $[m_{\sigma (1)},\dots ,m_{\sigma (n)}] = \Sign (\sigma )[m]$,
for any permutation $\sigma $.
\item $[m] = 0$ if $\det m = 0$.
\item $\sum _{i=0}^{n} (-1)^{i}[m_{0},\dots ,\hat m_{i},\dots ,m_{n}]
= 0$, for any $n+1$ vectors $m_{0},\dots ,m_{n}$ (the ``cocycle relation'').
\end{enumerate}

By the first relation, $\MMM _{n} $ is generated by those $[m]$ such
that $m_{i}$ is integral and primitive for all $i$.  If $m \in \SL_{n}
(\Z ) $, then $[m]$ is called a \emph{unimodular symbol}.  We have the
following fundamental result of Manin ($n=2$) and Ash-Rudolph ($n\geq
2$):

\begin{theorem}
\cite{manin,ash.rudolph} Any modular symbol is homologous to a finite
sum of unimodular symbols.
\end{theorem}

We sketch the proof.  If $|\det m |>1$, then one can show there
exists $v\in \Z ^{n}\smallsetminus \{0 \}$ such that 
\begin{equation}\label{decrease}
0 \leq  |\det m_{i} (v)| < |\det m|, \quad \hbox{for $i=1,\dots ,n$.}
\end{equation}
where $m_{i} (v)$ is the matrix obtained by replacing the column
$m_{i}$ with $v$.  Such a $v$ is called a \emph{reducing point} for
$m$.  Then applying the cocycle relation to the tuple $v,m_{1},\dots
,m_{n}$ yields an expression for $[m]$ in terms of the symbols $[m_{i}
(v)]$.  By induction this completes the proof.

This process of rewriting a modular symbol as a sum of unimodular
symbols is called the \emph{modular symbol algorithm}.  Using this
algorithm one can compute the action of the Hecke operators on
$H^{\nu} (\Gamma )$ as follows.  There are only finitely many
unimodular symbols mod $\Gamma $, and from them one can select a
subset dual to a basis of $H^{\nu} (\Gamma )$.  A Hecke operator acts
on the modular symbols by taking a unimodular symbol into a sum of
nonunimodular symbols.  Hence the modular symbol algorithm allows one
to compute the Hecke action on a basis, from which one can easily
compute the eigenvalues.

\section{Symplectic modular symbols}\label{sympms.section}
\subsection{}
For the first generalization we replace the linear group with the
symplectic group \cite{sympms}.  Let $V$ be a $2n$-dimensional $\Q
$-vector space with basis $\{e_{1},\dots ,e_{n},e_{\bar n},\dots
,e_{\bar 1} \}$, where $\bar \imath := 2n+1-i$.  Let
$\langle \phantom{a},\phantom{a}\rangle \colon V\times V\rightarrow \Q
$ be the nondegenerate, alternating bilinear form defined by 
\[
\langle e_{i}, e_{j} \rangle= \cases    {
    \begin{array}{ll} 1   & \quad \mbox{if $j = \bar\imath$ with $i<j$} \\
                      -1 & \quad\mbox{if $j = \bar\imath$ with $i>j$}\\
                       0& \quad\mbox{otherwise.}
    \end{array} }
\]
The form $\langle
\phantom{a},\phantom{a}\rangle $ is called a \emph{symplectic form},
and the \emph{symplectic group} $\Sp_{2n} (\Q )$ is defined to be the
subgroup of $\SL_{n} (\Q )$ preserving $\langle
\phantom{a},\phantom{a}\rangle $.

\subsection{}
Much of \S\ref{class.mod.sym} carries over without change, but there
are some new wrinkles coming from the geometry of the symplectic form.
Recall that an \emph{isotropic subspace} is one on which the
symplectic form vanishes, and that maximal (necessarily
$n$-dimensional) isotropic subspaces are called \emph{Lagrangian}.
Then the symplectic building $\BBB_{\Sp}$ has a $k$-simplex for every
length $(k+1)$ flag of isotropic subspaces.  Since the columns of a
symplectic matrix $m$ satisfy
\begin{equation}\label{isotropy}
\langle m_{i},m_{j}\rangle = 0 \quad \hbox{if and only if} \quad i\not
=\bar\jmath , 
\end{equation}
it is easy to see that $m$ determines $2^{n}\cdot n!$ oriented simplices of
maximal dimension in $\BBB_{\Sp}$.  

Furthermore, the arrangement of these simplices in $\BBB _{\Sp}$
differs from the linear case.  Suppose we use the columns of $m$ to
induce points in the projective space $\Proj ^{2n-1} (\Q )$.  Then the
Lagrangian subspaces spanned by the columns of $m$ become
$(n-1)$-dimensional flats arranged in the configuration of a
\emph{hyperoctahedron}.\footnote{Recall that a hyperoctahedron is the
convex hull of the $2n$ points $\{\pm e \mid e \in E\}$, where $E$ is
the standard basis of $\R ^{2n}$.}  This time $m$ determines a class
$[m]\in H_{n-1}  (\BBB_{\Sp})$, and as $m$ ranges over
all rational matrices with columns satisfying \eqref{isotropy}, the
duals of the classes $[m]$ span $H^{\nu} (\Gamma )$.

\subsection{}\label{cocycle.relation}
As a first step towards a symplectic modular symbol algorithm, one
must understand the analogues of the relations from \S\ref{rels}.  The
analogues of 1--3 are only slightly different to reflect the
hyperoctahedral symmetry.  
The cocycle relation, however, is more interesting.  A symbol $[m]$
and a generic nonzero rational point $v\in V$ determine $2n$ modular
symbols $[m_{i} (v)]$ as follows.  For any pair $(i,j)$ with $i\not =
\bar\jmath $, we define points $m_{ij}$ by
\[
m_{ij} := \langle v,m_{j}\rangle m_{i} - \langle v,m_{i}\rangle m_{j}.
\]
Let $[m_{i} (v)]$ be the modular symbol obtained by replacing $m_{\bar
\imath}$ with $v$, and replacing the $m_{j}$ with $j\not \in
\{i,\bar\imath \}$ by $m_{ij}$.  Then one can show $[m] = \sum
\varepsilon _{i}[m_{i} (v)]$ for appropriate signs $\varepsilon _{i}$.

For an example of this relation, consider Figure \ref{square.fig}.
The figure on the left shows the cocycle relation for $\Sp_{4}$ in terms of
a configuration in $\Proj ^{3}$.  The black dots are the points
corresponding to the $m_{i}$, the grey dot correspond to $v$, and the
triangles to the points $m_{ij}$.

\subsection{}
Now we can describe the symplectic modular symbol algorithm.  Let
$m\in M_{2n} (\Z )$ have columns satisfying \eqref{isotropy}.  Then
$\det m = \prod _{i=1}^{n}\langle m_{i}, m_{\bar \imath }\rangle$, and
one can show that if $|\det m| > 1$, there exists a vector $v\in \Z
^{n}\smallsetminus \{0 \}$ such that
\[
0 \leq |\langle m_{i}, v\rangle| < \langle m_{i}, m_{\bar \imath
}\rangle, \quad \hbox{for $i=1,\dots ,2n$.}
\] 

\begin{figure}[ht]
\centerline{\includegraphics[scale = .3]{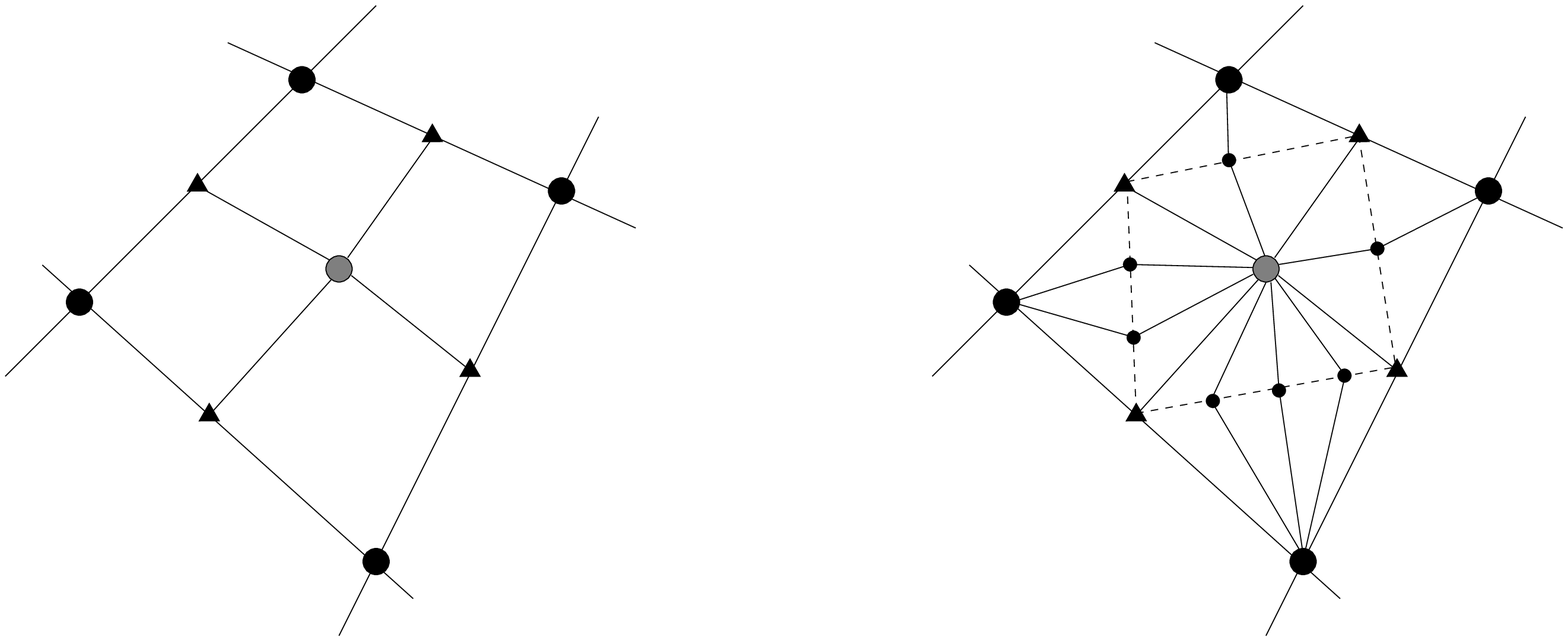}}
\caption{$G=\Sp _{4}$.  On the left, the outer square is the original
symbol $[m]$, and the four smaller squares are the symbols $[m_{i}
(v)]$.  On the right, each modular symbol has been further reduced by
applying the modular symbol algorithm to $\Sp _{2} = \SL _{2}$ modular
symbols.\label{square.fig}}
\end{figure}

We can apply $v$ to $[m]$ in the cocycle relation alluded to in
\S\ref{cocycle.relation}, but we will unfortunately find that $|\det
m_{i} (v)| > |\det m|$ in general.  However, all is not lost.  It
turns out that for fixed $i$ and fixed $v$, the $2n-2$ vectors
$\{m_{ij} \mid j\not= i, \bar\imath\}$ form a tuple that can be
regarded as a symplectic modular symbol associated to $\Sp_{2n-2}$.
By induction one knows how to make these symbols unimodular, and this
allows one to further reduce the $[m_{i} (v)]$ (cf. the right of
Figure \ref{square.fig}).

\section{Below the cohomological dimension}\label{below.coho}
\subsection{}
We return to the case of $\SL_{n}$.  As said before, a limitation of
the modular symbol algorithm is that one can compute the Hecke action
only on the top degree cohomology.  For $n\leq 3$ this cohomology
group is very interesting: it contains \emph{cuspidal} classes,
i.e. classes associated to cuspidal automorphic forms.  If $n\geq 4$,
however, the top degree cohomology group no longer contains cuspidal
classes.  In particular, if $n=4$, one is really interested in
computing the Hecke action on $H^{5} (\Gamma )$, and the modular
symbol algorithm applies to $H^{6} (\Gamma )$.

In this section we describe an algorithm that for $n\leq 4$ allows
computation of the Hecke action on $H^{\nu - 1} (\Gamma )$
\cite{experimental}.  However, there is one caveat: we cannot prove
the algorithm will terminate.  In practice, happily, the algorithm
has always converged, and has permitted investigation of this cohomology
\cite{computation}.

\subsection{}\label{questions}
To compute with lower degree cohomology groups, we use the
\emph{sharbly complex} $S_{*}$ \cite{ash.sharb}.  For $k\geq 0$, let
$S_{k}$ be the $\Z\Gamma $-module generated by the symbols
$\uu = [v_{1},\dots ,v_{n+k}]$, where $v_{i}\in \Q \smallsetminus \{0 \}$,
modulo the analogues of relations 1--3 in \S\ref{rels}.  Elements of
$S_{k}$ are called \emph{$k$-sharblies}.  Let $\partial \colon
S_{k}\rightarrow S_{k-1}$ be the map $\uu \mapsto
\sum _{i} (-1)^{i}[v_{1},\dots ,\widehat{v_{i}},\dots ,v_{n+k}]$,
linearly extended to all of $S_{k}$.  There is a map $S_{0}\rightarrow
\MMM _{n}$ giving a $\Z \Gamma $-free resolution of $\MMM _{n}$, and
one can show that this implies $H^{\nu -k} (\Gamma ; \C ) \cong H_{k}
(S_{*}\otimes \C )$.

As in \S\ref{rels}, it suffices to consider $k$-sharblies $\uu =
[v_{1},\dots ,v_{n+k}]$ with all $v_{i}$ integral and primitive.  Any
modular symbol of the form $[v_{i_{1}},\dots ,v_{i_{n}}]$, where
$\{i_{1},\dots ,i_{n} \}\subset \{1,\dots ,n+k \}$, is called a
\emph{submodular symbol} of $\uu $.  

Let $\xi =\sum n (\uu )\uu $ be a
sharbly chain.  We denote by $\| \xi \|$ the maximum absolute value of
the determinant of any submodular symbol of $\xi $.  The chain $\xi $
is called \emph{reduced} if $\| \xi \|=1$.  It is known
that reduced $1$-sharbly cycles provide a finite spanning set of
$H^{\nu-1 } (\Gamma;\C )$ for $n\leq 4$.

Since the Hecke operators take reduced sharbly cycles to nonreduced
cycles, our goal is to apply the modular symbol algorithm
\emph{simultaneously} over a nonreduced $1$-sharbly cycle $\xi $ to
lower the determinants of the submodular symbols.  Hence we are faced
with two problems: first, how do we combine reducing points with the
original $1$-sharbly $\xi $ to produce a new $1$-sharbly $\xi '$
homologous to $\xi $; second, how do we choose the reducing
points so that $\|\xi '\| < \|\xi \| $?

\subsection{}\label{construction}
To address the first issue we do the following.  Suppose $\uu =
[v_{1},\dots ,v_{n+1}]$ satisfies $n (\uu )\not =0$, and for
$i=1,\dots ,n+1$, let $\vv_{i} $ be the submodular symbol
$[v_{1},\dots ,\widehat{v_{i}},\dots ,v_{n+1}]$.  Assume that all
these submodular symbols are nonunimodular, and for each $i$ let $w_{i}$ be a
reducing point for $\vv_{i}$.  

For any subset $I\subset\{1,\dots ,n+1 \}$, let $\uu_{I}$ be the
$1$-sharbly $[u_{1},\dots ,u_{n+1}]$, where $u_{i}=w_{i}$ if $i\in I$,
and $u_{i}=v_{i}$ otherwise.  Then we have a relation in $S_{1}$ given
by
\begin{equation}\label{relation}
\uu = -\sum _{I\not =\varnothing} (-1)^{\#I}\uu_{I}.
\end{equation}
Geometrically this relation can be expressed using the combinatorics
of the hyperoctahedron \cite[\S4.4]{experimental}.  More generally, if
some $\vv_{i}$ happen to be unimodular, then one can construct a similar
relation using an iterated cone on a hyperoctahedron.

\subsection{}
Now we apply the construction in \S\ref{construction} to all the
$1$-sharblies $\uu $ with $n (\uu )\not =0$, and we choose reducing
points $\Gamma $-equivariantly.  Specifically, if $\vv$ and $\vv'$ are
two submodular symbols of $\xi $ with $\gamma \vv = \vv'$, then we choose
the corresponding reducing points such that $\gamma w = w'$.  After
applying \eqref{relation} to all the $\uu$ we determine a new
$1$-sharbly cycle $\xi '$.  Clearly $\xi '$ is homologous to $\xi $.
We claim that $\|\xi '\|$ should be less than $\|\xi \|$.

To see why this should be true, consider the $1$-sharblies $\uu_{I}$
on the right of \eqref{relation}.  Of these $1$-sharblies, those with
$\#I = 1$ contain the $\vv_{i}$ among their submodular symbols.  We
claim that since $\xi $ is a \emph{cycle} mod $\Gamma $, and since the
reducing points were chosen $\Gamma $-equivariantly over $\xi $, these
$1$-sharblies will not appear in $\xi '$.  Hence by construction we
have eliminated some of the ``bad'' submodular symbols from $\xi $.

\subsection{}
Unfortunately, this doesn't guarantee that $\|\xi '\|< \|\xi \|$.  The
problem is that we have no way of knowing that the submodular symbols
of the $\uu_{I}$ with $\#I>1$ don't have large determinants.  Indeed,
this brings us back to the second question raised in
\S\ref{questions}, since if the reducing points are chosen na\"\i
vely, these submodular symbols \emph{will} have large determinants.
However, we claim that one can (conjecturally) choose the reducing
points ``uniformly'' over $\xi $ in a sense by using LLL-reduction,
and that this problem doesn't occur in practice.  In fact, in
thousands of computer tests and in applications, we have always found
$\|\xi '\|< \|\xi \|$ if $n\leq 4$ and $\|\xi \|>1$.  We refer the
interested reader to \cite{experimental} for details.

\section{Self-adjoint homogeneous cones}\label{sahc.section}
\subsection{}
Now we describe a different approach to computing the Hecke action
that can be found in \cite{msa,sahc}.  The main idea is to replace
modular symbols and sharbly chains with chains built from rational
polyhedral cones, and to replace ``unimodularization'' with moving the
support of a chain into a certain canonically defined set of rational
polyhedral cones.  The results of this section apply to any arithmetic
group that is associated to a \emph{self-adjoint homogeneous cone};
the reduction theory in this generality is due to Ash
\cite[Ch. 2]{amrt}.  However, for simplicity we describe the results
in the context of \Vor 's work reduction theory of real
positive-definite quadratic forms \cite{voronoi1}.

Let $V$ be the real vector space of all real symmetric $n\times n$
matrices, and let $C$ be the subset of positive-definite matrices.
Then $C$ is a cone, i.e. $C$ is a convex set closed under
homotheties and containing no straight line.  The group $\SL_{n} (\Z
)$ acts on $V$ preserving $C$, and the action commutes with
homotheties.  In fact, modulo homotheties $C$ is isomorphic to $X =
\SL_{n} (\R )/\SO (n)$; this exhibits a hidden linear structure of the
symmetric space $X$.

Let $\bar C$ be the closure of $C$ in $V$.  \Vor\ showed how to
a set $\V $ of rational
polyhedral cones in $\bar C$ such that 
\begin{enumerate}
\item $\Gamma $ acts on $\V $.
\item If $\sigma \in \V $ then so is any face
of $\sigma $.
\item If $\sigma ,\tau \in \V $, then $\sigma \cap
\tau $ is a face of each.
\item Modulo $\Gamma $, the set $\V $ is finite.
\item The intersections $\sigma \cap C$ cover $C$. 
\end{enumerate}
The cones $\V $ provide a reduction theory for $C$ in the following
sense: any $x\in C$ lies in a unique cone $\sigma (x)\in \V $, and the
number of $\gamma \in \Gamma $ such that $\gamma \cdot \sigma (x) =
\sigma (x)$ is bounded.  Given $x\in C$, there is an explicit algorithm, the
\emph{\Vor \ reduction algorithm}, to find $\sigma (x)$.

The \Vor \ cones descend modulo homotheties to induce a decomposition
of $X$ into cells.  Furthermore, we can enlarge $C$ to a cone
$\tilde{C}$ such that, if $\tilde{X}$ denotes $\tilde{C}$ modulo
homotheties, then the quotient $\Gamma \backslash \tilde{X}$ is
compact.  This \emph{Satake compactification} of $\Gamma \backslash X$
is singular in general, but nevertheless can still be used to
compute $H^{*} (\Gamma ; \C )$.  For us, the salient points are that
the images of the \Vor \ cones induce a decomposition of $\tilde{C}$,
with all the properties listed above, and that the \Vor \ reduction
algorithm extends to the boundary $\partial \tilde{C} :=
\tilde{C}\smallsetminus C$.

\subsection{}
Now let $\cones ^{R}_{*}$ be the $\C $-complex generated by \emph{all}
simplicial rational polyhedral cones in $\tilde{C}$, and let $\cones
^{V}_{*}$ be the subcomplex generated by \Vor \
cones.\footnote{Although the \Vor \ cones aren't necessarily
simplicial, we can assume that they have been $\Gamma $-equivariantly
subdivided.}  For any chain $\xi \in \cones _{*}^{R}$, let $\supp \xi
$ be the set of cones supporting $\xi $.  The complex $\cones
^{R}_{*}$ is analogous to the sharbly complex, and the subcomplex
$\cones ^{V}_{*}$ to the subcomplex generated by the reduced
sharblies.  In general, however, $\cones ^{V}_{*}$ is not isomorphic
to the complex of reduced sharblies.  Cycles $\xi \in \cones ^{V}_{*}$
can be used to compute $H^{*}(\Gamma) $, but the image $T (\xi )$ of
$\xi $ under a Hecke operator won't be supported on \Vor \ cones.
Hence we must show how to push $T (\xi )$ back into $\cones ^{V}_{*}$.

To accomplish this we have essentially two tools---we can subdivide
the cones in $\supp T (\xi )$, and we can use the \Vor \ reduction
algorithm to determine the cone any point lies in.  We apply these as
follows.  Using the linear structure on $\tilde{C}$, we first
subdivide $T (\xi )$ very finely into a chain $\xi '$.  Then to each
$1$-cone $\tau \in \supp \xi '$, we assign a $1$-cone $\rho _{\tau
}\in \partial \tilde{C}$, and we use the combinatorics of $\xi '$ to
assemble the $\rho _{\tau }$ into a cycle $\xi ''$ homologous to $\xi
$.  We claim that if $\xi '$ is constructed so that $1$-cones
$\tau \in \supp \xi '$ lie in the \emph{same or adjacent \Vor \ cones,} then
the $\rho _{\tau }$ can be chosen to ensure $\xi '' \in \cones
^{V}_{*}$.

\subsection{}
We illustrate this process for $\SL_{2}$; more details can be found in
\cite{msa}.  Modulo homotheties the three-dimensional cone $\tilde{C}$
becomes the extended upper halfplane $\hp^{*}:= \hp \cup \Q \cup
\{\infty \} $, with $\partial \tilde{C}$ passing to the cusps $\hp
^{*}\smallsetminus \hp $.  The $3$-cones in $\V$ tiling $C$ pass to
the $\SL _{2} (\Z )$-translates of the ideal triangle with vertices at
$0,1,\infty $.  Let us call these ideal triangles \emph{\Vor \
triangles}.

\begin{figure}[ht]
\psfrag{mu}{$\mu $}
\psfrag{u1}{$u_{1}$}
\psfrag{u2}{$u_{2}$}
\centerline{\includegraphics[scale = .5]{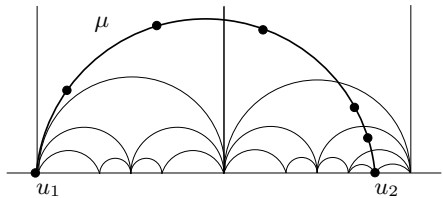}}
\caption{A subdivision of $\mu $; the solid dots are the $e_{i}$.
Since the $e_{i}$ lie in the same or adjacent \Vor \ triangles, we can
assign cusps to them to construct a homology to a cycle in $\cones ^{V}_{*}$.\label{partition.fig}}
\end{figure}

If $\xi \in \cones ^{R}_{*}$ is dual to a class in $H^{1} (\Gamma )$
and is supported on one $2$-cone, then $\supp \xi $ passes to a geodesic $\mu
$ between two cusps $u_{1}$, $u_{2}$ (Figure \ref{partition.fig}).  We
can subdivide $\mu $ into geodesic segments $\{\mu _{i} \}$ so that
the endpoints $e_{i}, e_{i+1}$ of $\mu _{i}$ lie in the same or
adjacent \Vor\ triangles.  Then we assign cusps to the $e_{i}$ as
follows.  If $e_{i}$ is not an endpoint of $\xi $, then we assign any
cusp $c_{i}$ of the \Vor \ triangle containing $e_{i}$.  Otherwise, if
$e_{i}=u_{1}$ or $u_{2}$ and hence is an endpoint of $\mu$, then we
assign $e_{i}$ to itself.  This determines a homology between $\xi $
and a chain $\xi ''$ supported on cones passing to the segments
$[c_{i}, c_{i+1}]$.  These cones are \Vor \ cones, and thus $\xi ''\in
\cones ^{V}_{*}$.

\section{Well-rounded retracts}\label{wrr.section}
\subsection{}
To conclude this article, we describe unpublished work of MacPherson
and McConnell \cite{bobnmark} that allows one to compute the Hecke
action on those $\Gamma $ for which a \emph{well-rounded retract} $W$
is available.  Again for simplicity we focus on $\Gamma \subset \SL
_{n} (\Z )$; our first task is to explain what $W$ is.

Let $V = \R ^{n}$ with the standard inner product preserved by $\SO
(n)$, and let $L\subset V$ be a lattice.  For any $v\in V$, write
$\|v\|$ for the length of $v$.  Let $m (L)$ be the minimal nonzero
length attained by any vector in $L$, and let $M (L) = \{v\in L\mid
\|v\| = m (L)\}$.  Then $L$ is said to be \emph{well-rounded} if $M
(L)$ spans $V$.

\subsection{}\label{retraction}
Consider the space of cosets $Y = \SL _{n} (\Z )\backslash \SL_{n} (\R
)$.  This space can be interpreted as the space of oriented lattices
in $\R ^{n}$ modulo homotheties.  Let $W\subset Y$ be the subset of
well-rounded lattices, and for any $j=0,\dots ,n$, let $Y_{j} = \{L\in
Y\mid \dim \Span M (L) \geq j \}$.  Clearly $Y_{0} = Y$ and $Y_{n}=W$.

According to Ash \cite{ash.wrr}, there is an $\SO (n)$-equivariant
retraction $r\colon Y\rightarrow W$ constructed as follows.  Let $L\in
Y_{j}$, and write $V=V_{1}\oplus V_{2}$, where $V_{1} = (\Span M
(L))\otimes \R$, and $V_{2}$ is the orthogonal complement of $V_{1}$.
For $0<\lambda \leq 1$, let $T (\lambda )$ be the linear
transformation $(v_{1},v_{2})\mapsto (v_{1}, \lambda v_{2})$, and let
$L[\lambda ]$ be the image of $L$ under $T (\lambda )$.  There is a
critical value $\lambda _{0}$ for which $\dim \Span M (L[\lambda ]) >
j$.  Then we can define $r_{j}\colon Y_{j}\rightarrow Y_{j+1}$ by
$r_{j} (L) = L[\lambda _{0}]$.  These retractions can be composed to
define the retraction $r\colon Y\rightarrow W$, and the space $W$ is
the well-rounded retract.

Since $r$ is $\SO (n)$-equivariant, it induces a retraction 
$\SL_{n} (\Z )\backslash \SL_{n} (\R )/\SO (n)\rightarrow W/\SO (n)$.
Moreover, $W$ can be given the structure of a locally-finite regular
cell-complex.  In a certain sense, these cells are dual to the \Vor \
cones from \S\ref{sahc.section}: \Vor \ cones of codimension $k$ are in
bijection with $W$-cells of dimension $k$.  The construction works if
$\Gamma $ is replaced with any finite-index subgroup of $\SL _{n} (\Z
)$, and hence one has a convenient topological model to study the
cohomology of any such $\Gamma $.

\subsection{}
Now we consider how the ideas used in the construction of $W$ can be
applied to compute the action of the Hecke operators on cohomology.
Let $d = (d_{1},\dots ,d_{n})$ be a tuple of strictly positive
integers, and let $g (d) \in \GL_{n} (\Q )$ be the diagonal matrix
with entries $d$.  Let $\Gamma ' := \Gamma \cap g^{-1}\Gamma g$.  The
\emph{Hecke correspondence} associated to this data is the diagram $(
c_{1},c_{2})\colon \Gamma' \backslash X \rightarrow \Gamma \backslash
X $, where the two maps are defined by $c_{1} (\Gamma 'x ) = \Gamma x$
and $c_{2} (\Gamma 'x) = \Gamma gx$.  In terms of the above
description, $c_{1}^{-1}\circ c_{2}$ is the (multivalued) map that
takes any lattice $L$ to the set of sublattices $\{M \subset L \mid
L/M \cong \Z /d_{1}\Z \oplus \cdots \oplus \Z /d_{n}\Z \}$.  A Hecke
correspondence induces a map $c_{1}^{*}\circ (c_{2})_*$ on cohomology
that is exactly a classical Hecke operator.  For example, if $n=2$,
$p$ is a prime, and $d = (1,p)$, then the induced Hecke operator is
the usual $T_{p}$.

\subsection{}
Fix a tuple $d$ and a pair of lattices $M\subset L$ as above.  Choose
$u\in [1,\infty )$.  For $v\in L$, let $\|\phantom{v}\|_{u}$ be
$\|v\|$ if $v\in M$, and $u\cdot \|v\|$ otherwise.  Now we can
consider the retraction $r$ described in \S\ref{retraction}, but using
$\|\phantom{v}\|_{u}$ instead of $\|\phantom{v}\|$ as the notion of
length.  When $u=1$, the result is the usual retract $W$.  But for $u
= u_{0}$ sufficiently large, only vectors in $M$ will be detected in
the retraction. Since $M$ is itself a lattice, we have $W_{u_{0}}
\cong W$.

These two complexes $W_{1}$ and $W_{u_{0}}$ appear in a larger complex
$\WWW$ that depends on $n$ and $d$ and is fibered over the interval
$[1,u_{0}]$ with fiber $W_{u}$.  The fibers $W_{1}$ and $W_{u_{0}}$
map to $W$ by the maps $c_{1}$ and $c_{2}$, respectively.  One
computes the action of the Hecke operator by lifting a class on
$\Gamma \backslash W$ to $\Gamma '\backslash \WWW$, pushing the lift
across $\Gamma '\backslash \WWW$ to the face $\Gamma \backslash
W_{u_{0}}$, and then pushing down via $c_{2}$ to $\Gamma  \backslash W$.    

\bibliographystyle{amsplain}
\bibliography{survey}

\providecommand{\bysame}{\leavevmode\hbox to3em{\hrulefill}\thinspace}
\begin{thebibliography}{10}

\bibitem{ash.wrr}
A.~Ash, \emph{Small-dimensional classifying spaces for arithmetic subgroups of
  general linear groups}, Duke Math. J. \textbf{51} (1984), 459--468.

\bibitem{ash.sharb}
\bysame, \emph{Unstable cohomology of ${SL}(n,\mathscr{O})$}, J. Algebra
  \textbf{167} (1994), no.~2, 330--342.

\bibitem{agg}
A.~Ash, D.~Grayson, and P.~Green, \emph{Computations of cuspidal cohomology of
  congruence subgroups of ${SL}_3(\mathbf{{Z}})$}, J. Number Theory \textbf{19}
  (1984), 412--436.

\bibitem{computation}
A.~Ash, P.~E. Gunnells, and M.~McConnell, \emph{Cohomology of congruence
  subgroups of ${S}{L}_4({\Z})$}, in preparation.

\bibitem{exp.ind}
A.~Ash and M.~McConnell, \emph{Experimental indications of three-dimensional
  galois representations from the cohomology of ${{SL}}(3,{\Z})$}, Experiment.
  Math. \textbf{1} (1992), no.~3, 209--223.

\bibitem{amrt}
A.~Ash, D.~Mumford, M.~Rapaport, and Y.~Tai., \emph{Smooth compactifications of
  locally symmetric varieties}, Math. Sci. Press, Brookline, Mass., 1975.

\bibitem{apt}
A.~Ash, R.~Pinch, and R.~Taylor, \emph{An $\widehat{A_4}$ extension of ${\Q}$
  attached to a non-selfdual automorphic form on ${GL}(3)$}, Math. Ann.
  \textbf{291} (1991), 753--766.

\bibitem{ash.rudolph}
A.~Ash and L.~Rudolph, \emph{The modular symbol and continued fractions in
  higher dimensions}, Invent. Math. \textbf{55} (1979), 241--250.

\bibitem{borel.serre}
A.~Borel and J.-P. Serre, \emph{Corners and arithmetic groups}, Comm. Math.
  Helv. \textbf{48} (1973), 436--491.

\bibitem{bygott}
J.~Bygott, \emph{Modular symbols and computation of cusp forms over imaginary
  quadratic fields}, Ph.D. thesis, Exeter University, 1997.

\bibitem{crem}
J.~E. Cremona, \emph{Hyperbolic tessellations, modular symbols, and elliptic
  curves over complex quadratic fields}, Compositio Math. \textbf{51} (1984),
  no.~3, 275--324.

\bibitem{msa}
P.~E. Gunnells, \emph{Modular symbols for ${\Q}$-rank one groups and
  {V}orono\u\i\ reduction}, J. Number Theory \textbf{75} (1999), no.~2,
  198--219.

\bibitem{experimental}
\bysame, \emph{Computing {H}ecke eigenvalues below the cohomological
  dimension}, Experiment. Math (to appear), 2000.

\bibitem{sympms}
\bysame, \emph{Symplectic modular symbols}, Duke Math. J., (to appear), 2000.

\bibitem{sahc}
P.~E. Gunnells and M.~McConnell, \emph{Hecke operators and ${\Q}$-groups
  associated to self-adjoint homogeneous cones}, math.NT/9811133, 1998.

\bibitem{bobnmark}
R.~MacPherson and M.~McConnell, \emph{Explicit reduction theory for {H}ecke
  correspondences}, in preparation.

\bibitem{manin}
Y.-I. Manin, \emph{Parabolic points and zeta-functions of modular curves},
  Math. USSR Izvestija \textbf{6} (1972), no.~1, 19--63.

\bibitem{merel}
L.~Merel, \emph{Universal {F}ourier expansions of modular forms}, On Artin's
  conjecture for odd $2$-dimensional representations, Springer, Berlin, 1994,
  pp.~59--94.

\bibitem{teit}
J.~T. Teitelbaum, \emph{Modular symbols for ${\F}\sb q({T})$}, Duke Math. J.
  \textbf{68} (1992), no.~2, 271--295.

\bibitem{geemen.top}
B.~{van Geemen} and J.~Top, \emph{A non-selfdual automorphic representation of
  ${GL} \sb 3$ and a {G}alois representation}, Invent. Math. \textbf{117}
  (1994), no.~3, 391--401.

\bibitem{vgt2}
B.~van Geemen, W.~van~der Kallen, J.~Top, and A.~Verberkmoes, \emph{Hecke
  eigenforms in the cohomology of congruence subgroups of
  $\rm{{S}{L}}(3,{{\Z}})$}, Experiment. Math. \textbf{6} (1997), no.~2,
  163--174.

\bibitem{voronoi1}
G.~Vorono\v\i, \emph{Sur quelques propri\'et\'es des formes quadratiques
  positives parfaites}, J. Reine Angew. Math. \textbf{133} (1908), 97--178.

\end{thebibliography}

\end{document}